\documentclass[11pt]{amsart}
\usepackage{amsmath,amsthm,amssymb,amscd}
\theoremstyle{plain}

\newtheorem{thm}{Theorem}[section]
\newtheorem{proposition}[thm]{Proposition}

\newtheorem{lem}[thm]{Lemma}

\theoremstyle{definition}

\newtheorem{defn}[thm]{Definition}
\newtheorem{ex}[thm]{Example}
\newtheorem{rem}[thm]{Remark}

\newcommand\cA{{\mathcal{A}}}

\newcommand\cD{{\mathcal{D}}}
\newcommand\cE{{\mathcal{E}}}
\newcommand\cF{{\mathcal{F}}}
\newcommand\cG{{\mathcal{G}}}
\newcommand\cH{{\mathcal{H}}}

\newcommand\cK{{\mathcal{K}}}
\newcommand\cL{{\mathcal{L}}}

\newcommand\cO{{\mathcal{O}}}

\newcommand\cR{{\mathcal{R}}}

\newcommand\RR{\mathbb{R}}

\newcommand\ZZ{\mathbb{Z}}

\newcommand\CC{\mathbb{C}}

\newcommand\End{\operatorname{End}}
\newcommand\supp{\mathop{\rm supp}\nolimits}
\newcommand{\cinf}{C^\infty}

\allowdisplaybreaks

\begin{document}
\title[The Egorov theorem for transverse Dirac type operators]{The Egorov theorem
for transverse Dirac type operators on foliated manifolds}
\author{Yuri A. Kordyukov}
\address{Institute of Mathematics, Russian Academy of Sciences, 112
Chernyshevsky street, 450077 Ufa, Russia}
\email{yurikor@matem.anrb.ru}

\thanks{Supported by the Russian Foundation of Basic Research
(grant no. 07-01-00081)}



\begin{abstract}
Egorov's theorem for trans\-ver\-sal\-ly elliptic operators, acting
on sections of a vector bundle over a compact foliated manifold, is
proved. This theorem relates the quantum evolution of transverse
pseudodifferential operators determined by a first order
transversally elliptic operator with the (classical) evolution of
its symbols determined by the parallel transport along the orbits of
the associated transverse bicharacteristic flow. For a particular
case of a transverse Dirac operator, the transverse bicharacteristic
flow is shown to be given by the transverse geodesic flow and the
parallel transport by the parallel transport determined by the
transverse Levi-Civita connection. These results allow us to
describe the noncommutative geodesic flow in noncommutative geometry
of Riemannian foliations.
\end{abstract}
\date{}
 \maketitle \hyphenation{trans-ver-sal-ly}

\bibliographystyle{plain}

\section*{Introduction}
The Egorov theorem is a fundamental fact in microlocal analysis and
quantum mechanics. It relates the evolution of pseudodifferential
operators on a compact manifold (quantum observables) determined by
a first order elliptic operator with the corresponding evolution of
classical observables --- the bicharacteristic flow on the space of
symbols. More precisely, let $M$ be a compact manifold and let $P$
be a positive, self-adjoint, elliptic, first order
pseudodifferential operator on $M$ with the positive principal
symbol $p\in S^1(T^*M\setminus 0)$. Let $f_t$ be the
bicharacteristic flow of the operator $P$, that is, the Hamiltonian
flow of $p$ on $T^*M$. Egorov's theorem \cite{egorov} states that,
for any pseudodifferential operator $A$ of order $0$ with the
principal symbol $a\in S^0(T^*M\setminus 0)$, the operator
$A(t)=e^{itP}Ae^{-itP}$ is a pseudodifferential operator of order
$0$. The principal symbol $a_t\in S^0(T^*M\setminus 0)$ of this
operator is given by the formula $$ a_t(x,\xi)=a(f_t(x,\xi)), \quad
(x,\xi)\in T^*M\setminus 0. $$ In the particular case
$P=\sqrt{\Delta_g}$, where $\Delta_g$ is the Laplace-Beltrami
operator of a Riemannian metric $g$ on $M$, the corresponding
bicharacteristic flow is the geodesic flow of $g$ on $T^*M$.

In \cite{Ja-Str06}, Egorov's theorem was extended to
pseudodifferential operators acting on sections of a vector bundle
$E$ on a compact manifold $M$. First, the authors gave an invariant
definition of the subprincipal symbol of a positive, self-adjoint,
elliptic, first order pseudodifferential operator $P\in \Psi^1(M,E)$
with the real scalar principal symbol $p\in S^1(T^*M\setminus 0)$ as
a partial connection along the Hamiltonian vector field of $p$ on
$T^*M$. The parallel transport along the orbits of the Hamiltonian
flow of $p$ on $T^*M$ defined by this partial connection determines
a flow $\beta_t$ acting on $S^0(T^*M\setminus 0, \End(\pi^*E))$.
Then the theorem in \cite{Ja-Str06} says that for any operator
$A\in\Psi^0(M,E)$ with the principal symbol $a\in S^0(T^*M\setminus
0, \End(\pi^*E))$, the operator $A(t)=e^{itP}Ae^{-itP}$ is in
$\Psi^0(M,E)$, and its principal symbol $a_t\in S^0(T^*M\setminus 0,
\End(\pi^*E))$ is given by $a_t=\beta_t(a)$.

If $P=\sqrt{\Delta_g}$, where $\Delta_g$ is the Hodge-Laplace
operator of a Riemannian metric $g$ acting on differential forms on
$M$, the corresponding flow $\beta_t$ on $S^0(T^*M\setminus 0,
\pi^*\End(\Lambda^*_{\CC}T^*M))$ is given by the parallel transport
along the orbits of the geodesic flow of $g$ on $T^*M$ with respect
to the Levi-Civita connection.

We also mention the works \cite{BG04,Dencker,Emmrich-W,GMMP,San99}
(and references therein) for discussion of Egorov's theorem for
matrix-valued operators and relations to parallel transport.

On the other side, in \cite{mpag} the author proved a version of
Egorov's theorem for scalar transversally elliptic operators on
compact foliated manifolds. For this purpose, we used the transverse
pseudodifferential calculus developed in \cite{noncom}. The
associated algebra of symbols is a noncommutative, Connes type
operator algebra associated with a natural foliation $\cF_N$ on the
conormal bundle $N^*\cF$ of the foliation $\cF$. The Egorov theorem
stated in \cite{mpag} relates the quantum evolution of transverse
pseudodifferential operators determined by a first order
transversally elliptic operator $P$ with the (classical) evolution
of its symbols determined by the transverse bicharacteristic flow of
$P$, which is the restriction of the bicharacteristic flow of $P$ to
$N^*\cF$. We also mention related works the Duistermaat-Guillemin
trace formula: \cite{dg:trace}  for transversally elliptic operators
on Riemannian foliations and \cite{San06} for the basic Laplacian of
a Riemannian foliation.

The main purpose of this paper is to extend Egorov's theorem to
transversally elliptic operators acting on sections of a holonomy
equivariant vector bundle on a compact foliated manifold, using
ideas of \cite{Ja-Str06}. In this case, it is shown that the
corresponding classical evolution is given by the parallel transport
along the orbits of the transverse bicharacteristic flow.
Furthermore, we introduce a natural class of first order
transversally elliptic operators, namely, transverse Dirac operators
$D_{\cE}$ with coefficients in an arbitrary holonomy equivariant
Hermitian vector bundle $\cE$, and compute the transverse
bicharacteristic flow for these operators. Quite remarkably, the
associated parallel transport is naturally determined by the
transverse Levi-Civita connection.

The transverse Dirac operators were introduced in \cite{G-K91}.
These papers mainly concern with the transverse Dirac operators
acting on basic sections (see also
\cite{G-K91a,G-K93,Jung01,Jung06,Jung07} and references therein).
The index theory of transverse Dirac operators was studied in
\cite{DGKY}. Finally, spectral triples defined by transverse
Dirac-type operators on Riemannian foliations were studied in
\cite{noncom,mpag}. In particular, the results of this paper can be
considered as a complement of our study of the noncommutative
geodesic flow started in \cite{mpag}.

I am grateful to D.Jakobson for sending to me a preprint version of
the paper \cite{Ja-Str06} before its publication and to the referee
for useful remarks.


\section{Classes of transverse pseudodifferential operators}\label{s:trpdo}
Throughout in the paper, $(M,{\mathcal F})$ is a compact foliated
manifold, $E$ is a Hermitian vector bundle on $M$,
$\operatorname{dim} M=n, \operatorname{dim} \cF=p, p+q=n$.

We will consider pseudodifferential operators, acting on
half-den\-si\-ti\-es. For any vector bundle $V$ on $M$, denote by
$|V|^{1/2}$ the associated half-density vector bundle. Let
$C^{\infty}(M,E)$ denote the space of smooth sections of the vector
bundle $E\otimes |TM|^{1/2}$, $L^2(M,E)$ the Hilbert space of square
integrable sections of $E\otimes |TM|^{1/2}$, ${\cD}'(M,E)$ the
space of distributional sections of $E\otimes |TM|^{1/2}$,
${\cD}'(M,E)=C^{\infty}(M,E)'$, and $H^s(M,E)$ the Sobolev space of
order $s$ of sections of $E\otimes |TM|^{1/2}$. Finally, let
$\Psi^{m}(M,E)$ denote the standard classes of pseudodifferential
operators, acting in $C^{\infty}(M,E)$.

We will use the classes $\Psi^{m,-\infty}(M,{\mathcal F},E)$ of
transversal pseudodifferential operators. Let us briefly recall its
definition, referring the reader to \cite{noncom} for more details.

First, for any $k_A \in S ^{m} (I ^{p} \times I^p\times I^q\times
{\RR}^{q}, {\cL}({\CC}^r))$ (here $r=\operatorname{rank} E$), define
an operator $A: C^\infty_c(I^n,\CC^r)\to C^\infty(I^n,\CC^r)$ by the
formula
\begin{equation}\label{loc}
Au(x,y)=(2\pi)^{-q} \int e^{i(y-y')\eta}k_A(x,x',y,\eta) u(x',y')
\,dx'\,dy'\,d\eta,
\end{equation}
where $u \in C^{\infty}_{c}(I^{n}, {\CC}^r), x \in I^{p}, y \in
I^{q}$. The function $k_A$ is called the complete symbol of $A$. As
usual, we will consider only classical (or polyhomogenous) symbols,
that is, those symbols, which can be represented as an asymptotic
sum of homogeneous (in $\eta$) components.

Let $\varkappa: U\subset M\rightarrow I^p\times I^q, \varkappa': U'
\subset M\rightarrow I^p\times I^q$ be a pair of compatible foliated
charts on $M$ equipped with trivializations of the bundle $E$ over
them. Any operator $A$ of the form~(\ref{loc}) with the Schwartz
kernel, compactly supported in $I^n\times I^n$, determines an
operator $A:C^{\infty}_c(U,\left. E\right|_U)\to
C^{\infty}_c(U',\left. E\right|_{U'})$, which extends to an operator
in $C^\infty(M,E)$ in a trivial way. The resulting operator is
called an elementary operator of class $\Psi
^{m,-\infty}(M,{\mathcal F},E)$.

\begin{defn}
The  class $\Psi ^{m,-\infty}(M,{\mathcal F},E)$ consists of all
operators $A$ in $C^{\infty}(M,E)$, which can be represented in the
form
\[
A=\sum_{i} A_i + K,
\]
where $A_i$ are elementary operators of class $\Psi
^{m,-\infty}(M,{\mathcal F},E)$, corresponding to a pair
$\varkappa_i,\varkappa'_i$ of compatible foliated charts, $K\in \Psi
^{-\infty}(M,E)$.
\end{defn}

\section{The principal symbol for transverse $\Psi$DOs}
In this Section, we will recall the definition of the principal
symbol for an operators of class $\Psi^{m,-\infty}(M,\cF,E)$.

First, we define the principal symbol of $A$ given by (\ref{loc}) as
the leafwise half-density
\begin{equation}\label{k-principal}
\sigma(A)(x,x^\prime,y,\eta)=k_{A,m}(x,x^\prime,y,\eta)|dx|^{1/2}|dx^\prime|^{1/2},
\end{equation}
where $k_{A,m}$ is the degree $m$ homogeneous component of the
complete symbol $k_A$.

Before giving the global definition of the principal symbol, we
recall several notions (for more details, see e.g. \cite{survey} and
references therein). Let $\gamma : [0.1]\to M$ be a continuous
leafwise path in $M$ with the initial point $x=\gamma(0)$ and the
final point $y=\gamma(1)$ and $T_0$ and $T_1$ arbitrary smooth
submanifolds (possibly, with boundary), transversal to the
foliation, such that $x\in T_0$ and $y\in T_1$. Sliding along the
leaves of the foliation $\cF$ determines a diffeomorphism
$H_{T_0T_1}(\gamma)$ of a neighborhood of $x$ in $T_0$ to a
neighborhood of $y$ in $T_1$, called the holonomy map along
$\gamma$. The differential of $H_{T_0T_1}(\gamma)$ at $x$ gives rise
to a well-defined linear map $T_xM/T_x\cF\to T_yM/T_y\cF$, which is
independent of the choice of transversals $T_0$ and $T_1$. This map
is called the linear holonomy map and denoted by
$dh_{\gamma}:T_xM/T_x\cF\to T_yM/T_y\cF$. The adjoint of
$dh_{\gamma}$ yields a linear map $dh_{\gamma}^*:N^*\cF_y\to
N^*\cF_x$, where we denote by $N^*{\mathcal F}$ the conormal bundle
to ${\mathcal F}$.

Denote by $G$ the holonomy groupoid of ${\cF}$. Recall that $G$
consists of $\sim_h$-equivalence classes of continuous leafwise
paths in $M$, where we set $\gamma_1\sim_h \gamma_2$, if $\gamma_1$
and $\gamma_2$ have the same initial and final points and the same
holonomy maps. $G$ is equipped with the source map $s:G\to M,
s(\gamma)=\gamma(0),$ and the range map $r:G\to M,
r(\gamma)=\gamma(1)$.

Let ${\mathcal F}_N$ be the linearized foliation in
$\tilde{N}^*{\mathcal F}=N^*\cF\setminus 0$ (cf., for instance,
\cite{Molino}). The leaf of the foliation $\cF_N$ through $\nu\in
\tilde{N}^*{\mathcal F}$ is the set of all points
$dh_{\gamma}^{*}(\nu)\in \tilde{N}^*{\mathcal F}$, where $\gamma\in
G, r(\gamma)=\pi(\nu)$ (here $\pi :T^*M\to M$ is the bundle map).

Consider a foliated chart $\varkappa: U\subset M\rightarrow
I^p\times I^q$ on $M$ with coordinates $(x,y)\in I^p\times I^q$ ($I$
is the open interval $(0,1)$) such that the restriction of $\cF$ to
$U$ is given by the sets $y={\rm const}$, equipped with a
trivialization of the vector bundle $E$. We will always assume that
the foliated chart $\varkappa$ is regular, that means that it admits
an extension to a foliated chart $\bar{\varkappa} : V\subset M
\stackrel{\sim}{\longrightarrow}(-2,2)^{n}$ with $\bar{U}\subset V$.
There is the corresponding chart in $T^*M$ with coordinates written
as $(x,y,\xi,\eta)\in I^p\times I^q\times \RR^p\times \RR^q$. In
these coordinates, the restriction of the conormal bundle
$N^*{\mathcal F}$ to $U$ is given by the equation $\xi=0$. So we
have a chart $\varkappa_n : U_1\subset N^*{\mathcal  F}
\stackrel{\sim}{\longrightarrow} I^p\times I^q\times \RR^q$ on
$N^*{\mathcal  F}$ with the coordinates $(x,y,\eta)\in I^p\times
I^q\times \RR^q$. This chart is a foliated chart on $N^*{\mathcal
F}$ for the linearized foliation ${\mathcal F}_N$, and the
restriction of ${\mathcal F}_N$ to $U_1$ is given by the level sets
$y= {\rm const}, \eta={\rm const}$.

The holonomy groupoid $G_{{\mathcal F}_N}$ of the linearized
foliation ${\mathcal F}_N$ can be described as the set of all
$(\gamma,\nu)\in G\times \tilde{N}^*{\mathcal F}$ such that
$r(\gamma)=\pi(\nu)$. The source map $s_N:G_{{\mathcal
F}_N}\rightarrow \tilde{N}^*{\mathcal F}$ and the range map
$r_N:G_{{\mathcal F}_N}\rightarrow \tilde{N}^*{\mathcal F}$ are
defined as $s_N(\gamma,\nu)=dh_{\gamma}^{*}(\nu)$ and
$r_N(\gamma,\nu)=\nu$. We have a map $\pi_G:G_{{\mathcal
F}_N}\rightarrow G$ given by $\pi_G(\gamma,\nu)=\gamma$.

The holonomy groupoid $G_{{\mathcal F}_N}$ carries a natural
codimension $q$ foliation $\cG_N$. The leaf of $\cG_N$ through a
point $(\gamma,\nu)\in G_{{\mathcal F}_N}$ is the set of all
$(\gamma', \nu')\in G_{{\mathcal F}_N}$ such that $\nu$ and $\nu'$
lie in the same leaf in $\cF_N$. Let $|T{\mathcal G}_N|^{1/2}$ be
the line bundle of leafwise half-densities on $G_{{\mathcal F}_N}$
with respect to the foliation ${\mathcal G}_N$. It is easy to see
that $$ |T{\mathcal G}_N|^{1/2}=r_N^*(|T{\mathcal
F}_N|^{1/2})\otimes s_N^*(|T{\mathcal F}_N|^{1/2}), $$ where
$s_N^*(|T{\mathcal F}_N|^{1/2})$ and $r_N^*(|T{\mathcal
F}_N|^{1/2})$ denote the lifts of the line bundle $|T{\mathcal
F}_N|^{1/2}$ of leafwise half-densities on $N^*{\mathcal F}$ via the
source and the range mappings $s_N$ and $r_N$ respectively.

Let $\pi^*E$ denote the lift of the vector bundle $E$ to
$\tilde{T}^*M=T^*M\setminus 0$ via the bundle map $\pi:\tilde{T}^*M
\to M$. Denote by ${\mathcal L}(\pi^*E)$ the vector bundle on
$G_{\cF_N}$, whose fiber at a point $(\gamma,\nu)\in G_{\cF_N}$ is
the space $\cL((\pi^*E)_{s_N(\gamma,\nu)},
(\pi^*E)_{r_N(\gamma,\nu)})$ of linear maps from
$(\pi^*E)_{s_N(\gamma,\nu)}$ to $(\pi^*E)_{r_N(\gamma,\nu)}$.

A section $k\in C^{\infty}(G_{{\mathcal F}_N}, {\mathcal
L}(\pi^*E)\otimes |T{\mathcal G}_N|^{1/2})$ is said to be properly
supported, if the restriction of the map $r:G_{\cF_N}\to
\tilde{N}^*\cF$ to $\supp k$ is a proper map. Consider the space
$C^{\infty}_{prop}(G_{{\mathcal F}_N}, {\mathcal L}(\pi^*E)\otimes
|T{\mathcal G}_N|^{1/2})$ of smooth, properly supported sections
of ${\mathcal L}(\pi^*E)\otimes |T{\mathcal G}_N|^{1/2}$. One can
introduce the structure of involutive algebra on
$C^{\infty}_{prop}(G_{{\mathcal F}_N}, {\mathcal L}(\pi^*E)\otimes
|T{\mathcal G}_N|^{1/2})$ by the standard formulas. Let
$S^{m}(G_{{\mathcal F}_N},{\mathcal L}(\pi^*E)\otimes |T{\mathcal
G}_N|^{1/2})$ be the space of all $s\in
C^{\infty}_{prop}(G_{{\mathcal F}_N},{\mathcal L}(\pi^*E)\otimes
|T{\mathcal G}_N|^{1/2})$ homogeneous of degree $m$ with respect
to the action of $\RR$ given by the multiplication in the fibers
of the vector bundle $\pi_G:G_{{\mathcal F}_N}\rightarrow G $.

Now let $\varkappa: U\subset M\rightarrow I^p\times I^q,
\varkappa': U' \subset M\rightarrow I^p\times I^q$, be two
compatible foliated charts on $M$. Then the corresponding foliated
charts $\varkappa_n: U_1\subset N^*{\mathcal  F}\rightarrow
I^p\times I^q \times \RR^q, \varkappa'_n: U'_1 \subset
N^*{\mathcal  F} \rightarrow I^p\times I^q \times \RR^q,$ are
compatible with respect to the foliation ${\mathcal  F}_N$. So
they define a foliated chart $V$ on the foliated manifold
$(G_{{\mathcal F}_N},{\mathcal G}_N)$ with the coordinates
$(x,x',y,\eta) \in I^p\times I^p\times I^q \times \RR^q$, and the
restriction of ${\mathcal  G}_N$ to $V$ is given by the level sets
$y= {\rm const}, \eta={\rm const}$.

The principal symbol $\sigma(A)$ of an operator $A\in
\Psi^{m,-\infty}(M,{\mathcal F},E)$ given in local coordinates by
the formula (\ref{k-principal}) is globally defined as an element
of the space $S^{m}(G_{{\mathcal F}_N},{\mathcal L}(\pi^*E)\otimes
|T{\mathcal G}_N|^{1/2})$. Thus, we have the half-density
principal symbol mapping
\begin{equation}\label{e:symbol}
\sigma: \Psi^{m,-\infty}(M,{\mathcal F},E)\rightarrow
S^m(G_{{\mathcal F}_N},{\mathcal L}(\pi^*E)\otimes |T{\mathcal
G}_N|^{1/2}),
\end{equation}
which satisfies
\[
\sigma(AB)=\sigma(A)\sigma(B),\quad \sigma(A^*)=\sigma(A)^*
\]
for any $A\in \Psi^{m_1,-\infty}(M,{\mathcal F},E)$ and $B\in
\Psi^{m_2,-\infty}(M,{\mathcal F},E)$.

\begin{ex}
Suppose that a foliation $\cF$ on a compact manifold $M$ is given by
the fibers of a fibration $f: M\to B$ over a compact manifold $B$.
Then, for any $x\in M$, $N^*_x\cF$ coincides with the image of the
cotangent map $f^* : T_{f(x)}^*B \to T_x^*M$. The inverse map
$(f^*)^{-1} : T_x^*M \to T_{f(x)}^*B$ is a fibration whose fibers
are the leaves of the linearized foliation $\cF_N$. Thus, we have
the diffeomorphism
\[
\{(x,\xi)\in M\times T^*B : f(x)=\pi_B(\xi)\} \stackrel{\cong}{\to}
N^*\cF, \quad (x,\xi)\mapsto f^*(\xi)\in N^*_x\cF,
\]
where $\pi_B :T^*B\to B$ is the cotangent bundle map. So the diagram
\[
  \begin{CD}
N^*\cF @>\pi>> M\\ @V(f^*)^{-1}VV @VVfV
\\ T^*B @>>\pi_B>
B
  \end{CD}
\]
commutes, and $N^*\cF$ can be considered as the pull-back of the
bundle $f: M\to B$ to $T^*B$:
\begin{equation}\label{e:nf}
N^*\cF\cong\pi_B^*(M)=\{ (x,\xi)\in M\times T^*B : f(x)=\pi(\xi) \}.
\end{equation}

The holonomy groupoid $G$ of $\cF$ is the fiber product
\[
M\times_BM=\{(x,y)\in M\times M : f(x)=f(y)\},
\]
where $s(x,y)=y, r(x,y)=x$. The holonomy groupoid $G_{\cF_N}$ can be
identified as above with
\[
N^*\cF\times_{T^*B}N^*\cF\cong \{(x,y,\xi)\in M \times M\times T^*B
: f(x)=f(y)=\pi_B(\xi)\},
\]
where $s_N(x,y,\xi)=y, r_N(x,y,\xi)=x$. So we have the commutative
diagram
\[
  \begin{CD}
G_{\cF_N} @>\pi_G >> G \\ @VVV @VVV
\\ T^*B @>>\pi_B>
B
  \end{CD}
\]
and the foliation $\cG_N$ is given by the fibers of the fibration
$G_{\cF_N}\to T^*B$.

For any $\xi\in T^*B$, let $\Psi^{-\infty}((N^*\cF)_\xi,
(\pi^*E)_\xi)$ be the involutive algebra of all smoothing operators,
acting on the space of smooth half-densities $C^\infty((N^*\cF)_\xi,
(\pi^*E)_\xi)$, where $(N^*\cF)_\xi$ is the fiber of the fibration
$N^*\cF\to T^*B$ at $\xi$ and $(\pi^*E)_\xi$ is the restriction of
$\pi^*E$ to $(N^*\cF)_\xi$. Consider a sheaf
$\Psi^{-\infty}(N^*\cF,\pi^*E)$ of involutive algebras on $T^*B$
whose stalk at $\xi\in T^*B$ is $\Psi^{-\infty}((N^*\cF)_\xi,
(\pi^*E)_\xi)$. For any section $\sigma$ of the sheaf
$\Psi^{-\infty}(N^*\cF,\pi^*E)$, the Schwartz kernels of the
operators $\sigma(\xi)$ determine a well-defined section of
${\mathcal L}(\pi^*E)\otimes |T{\mathcal G}_N|^{1/2}$ over
$G_{{\mathcal F}_N}\cong N^*\cF \times_{T^*B} N^*\cF$. We say that
$\sigma$ is smooth, if the corresponding section is smooth. This
defines an algebra isomorphism of $\Psi^{-\infty}(N^*\cF,\pi^*E)$
with $C^{\infty}(G_{{\mathcal F}_N}, {\mathcal L}(\pi^*E)\otimes
|T{\mathcal G}_N|^{1/2})$.
\end{ex}

\begin{rem}
Suppose as above that a foliation $\cF$ on a compact manifold $M$ is
given by the fibers of a fibration $f: M\to B$ over a compact
manifold $B$. If we consider the cotangent bundle $T^*M$ as a
symplectic manifold equipped with the canonical symplectic
structure, then $N^*\cF$ is a closed coisotropic submanifold, and
the linearized foliation $\cF_N$ coincides with the null-foliation
of this coisotropic submanifold, that is, $T{\cF}_N$ is the
skew-orthogonal complement of $T(N^*\cF)$ in $T(T^*M)$. It is
well-known that the fiber product $N^*\cF\times_{T^*B}N^*\cF$ is a
canonical relation in ${T}^*M$, which is often called the flowout of
the coisotropic submanifold ${N}^*{\cF}$. The algebra of Fourier
integral operators associated with this canonical relation was
introduced in \cite{GS79}. In this particular case, it coincides
with the algebra $\Psi^{*,-\infty}(M,{\cF},E)$.

For an arbitrary compact foliated manifold $(M,\cF)$, one can
consider $G_{{\cF}_N}$ as an immersed canonical relation in
${T}^*M$, and the associated algebra of Fourier integral operators
also coincides with $\Psi^{*,-\infty}(M,{\cF},E)$. One has to be
only a little bit careful, defining the algebra of Fourier integral
operators associated with an immersed canonical relation (see
\cite{noncom} for more details).
\end{rem}

\section{Transverse principal and subprincipal symbols}\label{s:pdo}
Recall that the principal symbol of an operator $P\in\Psi^m(M,E)$ is
an element of the space $S^m(\tilde{T}^*M, \End(\pi^*E))$ of smooth
sections of the vector bundle $\End(\pi^*E)$, homogeneous of degree
$m$ with respect to the $\RR$-multiplication in the fibers of
$\End(\pi^*E)$.

By definition, the transversal principal symbol $\sigma(P)$ of
$P\in\Psi^m(M,E)$ is the restriction of its principal symbol to
$\tilde{N}^*{\mathcal F}$. So we have
\[
\sigma(P)\in S^m(\tilde{N}^*\cF, \End(\pi^*E)).
\]

The principal symbol of $P$ in a foliated chart is given by the top
degree homogeneous component $p_m$ of its complete symbol $p$, and
the transverse principal symbol is given by
\[
\sigma(P)(x,y,\eta)=p_m(x,y,0,\eta), \quad (x,y,\eta)\in I^p\times
I^q\times \RR^q.
\]

Before passing to the definition of the transverse subprincipal
symbol, we recall the concept of a partial connection.

By a partial connection on a vector bundle $V$ over a smooth
manifold $X$ along a vector field $v$ on $X$ we will understand a
linear map $\nabla_v:C^\infty(X,V)\to C^\infty(X,V)$ satisfying
\[
\nabla_v(fs)=v(f)s+f\nabla_v(s), \quad f\in C^\infty(X),\quad s\in
C^\infty(X,V).
\]
If we fix a trivialization of $V$ over an open subset $U\subset M$,
then one can write
\[
\nabla_v = v \cdot Id + \Gamma
\]
on $C^\infty(U,\CC^N)$ for some $\Gamma\in C^\infty(U,
\End(\CC^N))$. Under a change of trivializations by a function $T\in
C^\infty(U, \operatorname{GL}(N,\CC))$, we get the transformation
law
\[
\Gamma'=T^{-1}\Gamma T+ T^{-1}v(T).
\]
Let $f_t : X\to X$ be the flow on $X$ generated by $v$. One can
define the parallel transport on $V$ along the orbits of $v$ as
follows. Let $x\in X$ and $w\in V_x$. Let the section $\tau \in
[0,t]\mapsto w(\tau)\in V_{f_\tau(x)}$ is a solution in local
coordinates of the Cauchy problem
\[
\frac{dw(\tau)}{d\tau}=\Gamma(f_\tau(x)), \quad \tau\in [0,t],
\]
\[
v(0)=w.
\]
The parallel transport of $w$ along the orbit $\{f_\tau(x) : \tau\in
[0,t]\}$ is defined as $\alpha_t(w)=w(t)\in V_{f_t(x)}$.

The parallel transport determines a flow $\alpha_t$ on the vector
bundle $V$ which projects to the flow $f_t$ on $X$ under the bundle
map $V\to X$ and makes $V$ an $\RR$-equivariant vector bundle.

The induced flow $\alpha^*_t$ on $C^\infty(X,V)$ satisfies
\begin{equation}\label{e:subflow}
\frac{d}{dt}\alpha^*_ts=\nabla_v(\alpha^*_ts), \quad s\in
C^\infty(X,V).
\end{equation}

Now we go back to the foliation setting. Assume that an operator
$P\in\Psi^m(M,E)$ has the scalar and real principal symbol $p_m\in
C^\infty(\tilde{T}^*M)$. Let $X_{p_m}$ be the Hamiltonian vector
field of $p_m$. Recall that, in a foliation chart, $X_{p_m}$ is
given by
\[
X_{p_m} = \sum_{j=1}^p \left(\partial_{\xi_j} p_{m} \frac{\partial
}{\partial x_j} - \partial_{x_j} p_{m} \frac{\partial }{\partial
\xi_j}\right) + \sum_{k=1}^q \left(\partial_{\eta_k}
p_{m}\frac{\partial }{\partial y_k} - \partial_{y_k} p_{m}
\frac{\partial }{\partial \eta_k}\right).
\]

As in \cite{Ja-Str06}, we define the subprincipal symbol of an
operator $P\in\Psi^m(M,E)$ as a partial connection $\nabla_{sub}(P)$
on $\pi^*E$ along the Hamiltonian vector field $X_{p_m}$. In local
coordinates, we have
\[
\nabla_{sub}(P)=X_{p_m}+ip_{sub},
\]
where
\begin{multline*}
p_{sub}(x,y,\xi,\eta) =p_{m-1}(x,y,\xi,\eta)\\ -
\frac{1}{2i}\sum_{j=1}^p\frac{\partial^2p_m}{\partial x_j\partial
\xi_j}(x,y,\xi,\eta)-
\frac{1}{2i}\sum_{l=1}^q\frac{\partial^2p_m}{\partial y_l\partial
\eta_l}(x,y,\xi,\eta).
\end{multline*}

Now in addition assume that the transverse principal symbol of
$P\in\Psi^m(M,E)$ is holonomy invariant. A function $\sigma\in
C^{\infty}(\tilde{N}^*{\mathcal F})$ is called holonomy invariant if
it satisfies the following condition:
\[
\sigma(dh^*_{\gamma}(\nu))=\sigma(\nu),\quad \gamma\in G, \quad
\nu\in N^*_{r(\gamma)}\cF.
\]
In a foliation chart, holonomy invariance of $\sigma$ means that
\[
\sigma(x,y,\eta)=\sigma(y,\eta), \quad (x,y,\eta)\in I^p\times
I^q\times \RR^q.
\]
Observe also that $\sigma\in C^{\infty}(\tilde{N}^*{\mathcal F})$ is
holonomy invariant if and only if it is constant along the leaves of
$\cF_N$.

Under these assumptions, the vector field $X_{p_m}$ is tangent to
$N^*\cF$. We define the transverse subprincipal symbol of $P$ as the
restriction of its subprincipal symbol to $\tilde{N}^*\cF$. In a
foliated coordinate chart, it is given by
\[
\nabla_{sub}(P)=X_{p_m}+i\sigma_{sub}(P),
\]
where
\begin{multline}\label{e:subprincipal1}
\sigma_{sub}(P)(x,y,\eta) =p_{m-1}(x,y,0,\eta)\\ -
\frac{1}{2i}\sum_{j=1}^p\frac{\partial^2p_m}{\partial x_j\partial
\xi_j}(x,y,0,\eta)-
\frac{1}{2i}\sum_{l=1}^q\frac{\partial^2\sigma(P)}{\partial
y_l\partial \eta_l}(y,\eta).
\end{multline}

\section{A $\Psi^{*}(M,E)$-bimodule structure}\label{s:subprincipal}
In this Section, we will study the structure of a
$\Psi^{*}(M,E)$-bimodule on the algebra
$\Psi^{*,-\infty}(M,{\mathcal F},E)$ given by the composition of
operators. An important new statement is the corresponding formula
for the complete symbols, which is given in the following theorem.

\begin{thm}\label{t:compose}
If $A$ is given by (\ref{loc}) with some $k_A\in S^{m_1} (I ^{p}
\times I^p\times I^q\times {\RR}^{q}, {\cL}({\CC}^r))$ and $B\in
\Psi^{m_2}(I^n, \CC^r)$, then $AB$ and $BA$ are given by (\ref{loc})
with some $k_{AB}\in S^{m_1+m_2} (I ^{p} \times I^p\times I^q\times
{\RR}^{q}, {\cL}({\CC}^r))$ and $k_{BA}\in S^{m_1+m_2} (I ^{p}
\times I^p\times I^q\times {\RR}^{q}, {\cL}({\CC}^r))$, which admit
the following asymptotic expansions
\begin{align*}
k_{AB}(x,x^\prime,y,\eta) & \sim
\sum_{\alpha,\beta}\frac{1}{\alpha ! \beta !} \partial_\xi^\alpha
\partial^\beta_\eta b(x,y,0,\eta)
 D^\alpha_{x} D^\beta_yk_{A}(x,x^\prime,y,\eta),\\
k_{BA}(x,x^\prime,y,\eta) & \sim
\sum_{\alpha,\beta}\frac{1}{\alpha ! \beta !}
D^\alpha_{x^\prime}\partial^\beta_\eta
k_{A}(x,x^\prime,y,\eta)(-\partial_\xi)^\alpha D^\beta_yb
(x^\prime,y,0,\eta).
\end{align*}
\end{thm}

The proof of this theorem can be achieved by a straightforward
modification of the standard arguments. As an immediate consequence,
we get:

\begin{proposition}[\cite{noncom}]\label{p:comm0}
\label{module} If $A \in  \Psi ^{m_{1},-\infty}(M,{\cF},E)$ and
$B\in \Psi ^{m_{2}}(M,E)$, then $AB$ and $BA$ in $\Psi
^{m_{1}+m_{2},-\infty}(M,{\cF},E)$ and
\[
\sigma(AB)=\sigma(A)\cdot r^*_N\sigma(B),\quad \sigma(BA)= s^*_N
\sigma(B)\cdot \sigma(A).
\]
\end{proposition}

Now we assume that $B\in \Psi^{m_2}(M, E)$ is such that the
principal symbol of is real and scalar, and its transverse principal
symbol is holonomy invariant. By Proposition~\ref{p:comm0}, it
follows that, for any $A\in \Psi^{m_1}(M, \cF, E)$, the operator
$[A,B]$ belongs to $\Psi ^{m_{1}+m_{2}-1,-\infty}(M,{\cF},E)$. Using
Theorem~\ref{t:compose}, one can compute  the principal symbol of
$[A,B]$.

Denote by $b_{m_2}$ the principal symbol of $B$. As above, $X_b$
denotes the restriction of the Hamiltonian vector field of $b_{m_2}$
to $N^*\cF$. Since $X_b$ is an infinitesimal transformation of
$\cF_N$, there exists a vector field $\cH_b$ on $G_{{\mathcal F}_N}$
such that $ds_N(\cH_b)=X_b$ and $dr_N(\cH_b)=X_b$. In local
coordinates, $\cH_b$ is given by
\begin{multline*}
\cH_b(x,x',y,\eta) = \sum_{j=1}^p \partial_{\xi_j}
b_{m_2}(x,y,0,\eta) \frac{\partial }{\partial x_j} +\sum_{j=1}^p
\partial_{\xi_j} b_{m_2}(x^\prime,y,0,\eta) \frac{\partial
}{\partial x^\prime_j}
\\ + \sum_{k=1}^q \left(\partial_{\eta_k} \sigma_B(y,\eta)\frac{\partial
}{\partial y_k} -
\partial_{y_k} \sigma_B(y,\eta) \frac{\partial }{\partial
\eta_k}\right),\\ \quad (x,x',y,\eta) \in I^p\times I^p\times
I^q\times \RR^q.
\end{multline*}

Denote by $\cL_{\cH_b}$ the Lie derivative on the space
$C^{\infty}_{prop}(G_{{\mathcal F}_N}, |T{\mathcal G}_N|^{1/2})$ by
the vector field $\cH_b$. In a foliated chart, it defines a
derivative on the space $C^{\infty}_{prop}(G_{{\mathcal F}_N},
{\mathcal L}(\pi^*E)\otimes |T{\mathcal G}_N|^{1/2})$. For any $k\in
C^{\infty}_{prop}(G_{{\mathcal F}_N}, {\mathcal L}(\pi^*E)\otimes
|T{\mathcal G}_N|^{1/2})$ of the form $k=k(x,x',y,\eta) |dx|^{1/2}
|dx^\prime|^{1/2}$, we have
\begin{multline*}
\cL_{\cH_b}k= \Big(\cH_bk(x,x',y,\eta)+\frac{1}{2}\sum_{j=1}^p
D_{x_j}\partial_{\xi_j} b_{m_2}(x,y,0,\eta) k(x,x^\prime,y,\eta) \\
-\frac{1}{2}\sum_{j=1}^p D_{x^\prime_j}\partial_{\xi_j}
b_{m_2}(x^\prime,y,0,\eta)
k(x,x^\prime,y,\eta)\Big) |dx|^{1/2}|dx^\prime|^{1/2},\\
\quad (x,x',y,\eta) \in I^p\times I^p\times I^q\times \RR^q.
\end{multline*}

The transverse subprincipal symbol of the operator $B$ considered as
a partial connection on $N^*\cF$ along $X_b$ yields the
corresponding partial connection on the space
$C^{\infty}_{prop}(G_{{\mathcal F}_N}, {\mathcal L}(\pi^*E)\otimes
|T{\mathcal G}_N|^{1/2})$ along $\cH_b$. In a foliation chart, for
any $k\in C^{\infty}_{prop}(G_{{\mathcal F}_N}, {\mathcal
L}(\pi^*E)\otimes |T{\mathcal G}_N|^{1/2})$, we have
\[
\nabla_{\cH_{b}}k=\cL_{\cH_b}k+i(k\cdot r^*_N\sigma_{sub}(B)-
s^*_N\sigma_{sub}(B)\cdot k).
\]

By a straightforward calculation, Theorem~\ref{t:compose} implies
the following result.

\begin{thm}\label{t:com}
Let $A\in \Psi^{m_1}(M, \cF, E)$ and $B\in \Psi^{m_2}(M, E)$.
Suppose that the principal symbol of $B$ is real and scalar, and the
transverse principal symbol of $B$ is holonomy invariant. Then
\begin{equation}\label{e:comm2}
    \sigma([B,A])=\frac{1}{i}\nabla_{\cH_{b}}\sigma(A).
\end{equation}
\end{thm}

\section{Transverse bicharacteristic flow}\label{s:flow}
In this Section, we give a definition of the transverse
bicharacteristic flow associated with a first order transversally
elliptic operator.

Consider an operator $P\in\Psi^1(M,E)$ which has the real scalar
principal symbol and the holonomy invariant transverse principal
symbol. Let $p\in S^1(\tilde{T}^*M)$ be the principal symbol of $P$.
The Hamiltonian flow $f_t$ of $p$ preserves $\tilde{N}^*{\mathcal
F}$, and its restriction to $N^*{\mathcal F}$ (denoted also by
$f_t$) preserves the foliation $\cF_N$, that is, takes any leaf of
$\cF_N$ to a leaf. Moreover, one can show that there exists a flow
$F_t$ on $G_{{\mathcal F}_N}$ such that $s_N\circ F_t=f_t\circ s_N$,
$r_N\circ F_t=f_t\circ r_N$, which preserves the foliation $\cG_N$.
Actually, this flow is generated by the vector field $\cH_p$
introduced in Section~\ref{s:subprincipal}. It is easy to see that
the flow $F_t$ depends only on the $1$-jet of the principal symbol
of $P$ along $N^*\cF$.

Let $\alpha_t^*$ be the flow on $C^\infty(N^*\cF,\pi^*E)$ determined
by the subprincipal symbol $\nabla^{sub}(P)$ of $P$ (see
(\ref{e:subflow})). It induces the flow
$\operatorname{Ad}(\alpha_t)^*$ on the space
$C^{\infty}_{prop}(G_{{\mathcal F}_N}, {\mathcal L}(\pi^*E)\otimes
|T{\mathcal G}_N|^{1/2})$, which satisfies
\[
\frac{d}{dt}\operatorname{Ad}(\alpha_t)^*k=\nabla_{\cH_{p}}k, \quad
k\in C^{\infty}_{prop}(G_{{\mathcal F}_N}, {\mathcal
L}(\pi^*E)\otimes |T{\mathcal G}_N|^{1/2}).
\]
This flow will be called the transverse bicharacteristic flow of
$P$. One can show that
\[
\operatorname{Ad}(\alpha_t)^*\circ s^*_N=s^*_N\circ \alpha_t^*,
\quad \operatorname{Ad}(\alpha_t)^*\circ r^*_N=r^*_N\circ
\alpha_t^*.
\]

Now consider a transversally elliptic operator $A\in\Psi^2(M,E)$,
which has the positive scalar principal symbol and the holonomy
invariant transverse principal symbol. (Recall that an operator
$P\in\Psi^m(M,E)$ is said to be transversally elliptic, if
$\sigma_P(\nu)$ is invertible for any $\nu\in\tilde{N}^*{\mathcal
F}$.) Let $a_2\in S^2(\tilde{T}^*M)$ be the principal symbol of $A$:
$a_2\geq 0$. Then the operator $\sqrt{A}$ is not, in general, well
defined, and even if $A$ is positive self-adjoint and the operator
$\sqrt{A}$ is a well defined positive operator in $L^2(M,E)$, it is
not, in general, a pseudodifferential operator. Nevertheless, we can
define its transverse bicharacteristic flow, working at the level of
symbols.

By assumption, $a_2$ is positive in some conic neighborhood of
$\tilde{N}^*{\mathcal F}$. Take any scalar elliptic symbol ${\tilde
p}\in S^1(\tilde{T}^*M)$, which is equal to $\sqrt{a_2}$ in some
conic neighborhood of $\tilde{N}^*{\mathcal F}$ (indeed, it is
sufficient that the 1-jets of ${\tilde p}$ and $\sqrt{a_2}$ coincide
on $\tilde{N}^*{\mathcal F}$). Proceeding as above, we obtain the
flow $\operatorname{Ad}(\alpha_t)^*$ on
$C^{\infty}_{prop}(G_{{\mathcal F}_N}, {\mathcal L}(\pi^*E)\otimes
|T{\mathcal G}_N|^{1/2})$, which is independent of a choice of
$\tilde{p}$ and will be called the transverse bicharacteristic flow
of $\sqrt{A}$.

\begin{ex}
Suppose that $\cF$ is a Riemannian foliation and $g_M$ is a
bundle-like metric on $M$. Recall that a Riemannian metric $g_M$ on
$M$ is bundle-like, if the induced metric on the normal bundle
$Q=TM/T\cF$ is holonomy invariant, that is, for any continuous
leafwise path $\gamma$ from $x$ to $y$, the corresponding linear
holonomy map $dh_\gamma : Q_x\to Q_y$ is an isometry (see, for
instance, \cite{Molino,Re} for more details on Riemannian
foliations).

For any $x\in M$, let $T^H_xM=T_x{\cF}^{\bot}$. So we have a smooth
vector subbundle $T^HM$ of $TM$ such that
\begin{equation}\label{e:decomp}
TM=T^HM\oplus T\cF.
\end{equation}
There is a natural isomorphism $T^HM\cong Q$. Observe also natural
isomorphisms $T^HM^*\cong Q^*\cong N^*\cF$.

The decomposition (\ref{e:decomp}) induces a bigrading on $\Lambda
T^{*}M$: $$ \Lambda^k
T^{*}M=\bigoplus_{i=0}^{k}\Lambda^{i,k-i}T^{*}M,\quad k=0,1,\ldots,
n,
$$ where $\Lambda^{i,j}T^{*}M=\Lambda^{i}T\cF^{*}\otimes \Lambda^{j}
T^HM^{*}$. In this bigrading, the de Rham differential $d$ can be
written as
$$ d=d_F+d_H+\theta, $$ where $d_F$ and $d_H$ are first order
differential operators (called the tangential de Rham differential
and the transversal de Rham differential accordingly), and $\theta$
is a zero order differential operator.

By definition, the transverse signature operator is a first order
differential operator in  $C^{\infty}(M,\Lambda T^HM^{*})$ given by
\[
D_H=d_H + d^*_H.
\]

The principal symbol of $D^2_H$ (see Theorem~\ref{t:trsym} below) is
given by
\[
a_2(x,\xi)=g^M(P^H(x,\xi), P^H(x,\xi)),\quad (x,\xi)\in T^*M,
\]
where $g^M$ is the induced metric on $T^*M$, $P^H: T^*M\to T^HM^*$
is the orthogonal projection. The holonomy invariance of the
transverse principal symbol is equivalent to the bundle-like
property of the metric.

The transverse bicharacteristic flow of the operator $\langle
D_H\rangle=(D^2+I)^{1/2}$ coincides with the transverse geodesic
flow $\gamma^M_t$ of $g_M$, which is the restriction of the geodesic
flow of $g_M$ to $N^*\cF$.
\end{ex}

\begin{ex}
Suppose that a foliation $\cF$ on a compact manifold $M$ is given by
the fibers of a fibration $f: M\to B$ over a compact manifold $B$. A
Riemannian metric $g_M$ on $M$ is bundle-like if and only if there
exists a Riemannian metric $g_B$ on $B$ such that, for any $x\in M$,
the tangent map $f_*$ induces an isometry from $(T^H_xM,
g_M|_{T^HM})$ to $(T_{f(x)}B, g_B)$, or, equivalently, $f :
(M,g_M)\to (B,g_B)$ is a Riemannian submersion. Then the transverse
geodesic flow $\gamma^M_t$ of $g_M$ projects under the map $f^*$ to
the geodesic flow $\gamma^B_t$ of $g_B$ that implies commutativity
of the following diagram
\[
  \begin{CD}
N^*\cF @>\gamma^M_t>> N^*\cF\\ @Af^*AA @AAf^*A
\\ T^*B @>\gamma^B_t>>
T^*B
  \end{CD}
\]
Commutativity of this diagram allows us to lift the flow
$\gamma^M_t$ to the holonomy groupoid $G_{\cF_N}\cong
N^*\cF\times_{T^*B}N^*\cF$ as above.
\end{ex}

\begin{ex}
Suppose that, in the setting of the previous example, the fibration
$f: M\to B$ is a principal $K$-bundle with a compact group $K$. The
group $K$ has a natural Hamiltonian action on the cotangent bundle
$T^*M$. The conormal bundle $N^*\cF$ is a $K$-invariant submanifold
of $T^*M$, and the fibration $(f^*)^{-1} : N^*\cF \to T^*B$ is a
principal $K$-bundle.

Suppose that $\omega$ is a connection on the principal bundle
$f:M\to B$. It gives rise to a decomposition
\[
T_mM=V_m\oplus H_m, \quad m\in M,
\]
where $V_m$ is the vertical space and $H_m$ is the connection's
horizontal distribution. The vertical space $V_m$ is naturally
isomorphic to the Lie algebra $\mathfrak{k}$ of $K$, and the
horizontal space $H_m$ is identified with the tangent space
$T_{f(m)}B$ to the base. Choose a Riemannian metric on $B$ and a
bi-invariant metric on $K$, and define a $K$-invariant Riemannian
metric on $M$, by requiring that, on $V_m$, it is induced by the
fixed bi-invariant metric on $K$, on $H_m$, it is the lift of the
Riemannian metric on $B$, and $V_m$ and $H_m$ are orthogonal. Such a
metric is sometimes called the Kaluza-Klein metric of the
connection. The fibers of the bundle $f:M\to B$ are totally geodesic
submanifolds, which are isometric to $K$.

The pull back of the connection form $\omega$ on $M$ defines a
connection form $\pi_M^*\omega$ on the principal bundle $(f^*)^{-1}:
N^*\cF\to T^*B$. The transverse geodesic flow $\gamma_t^M$ of the
Kaluza-Klein metric is described as follows. For any $\nu\in
N_m^*\cF,$ the element $\gamma_t^M(\nu)\in N^*\cF$ is obtained by
the parallel transport of $\nu$ along the orbit
$\{\gamma^B_\tau((f^*)^{-1}(\nu)) : \tau \in [0,t]\}$ of the
geodesic flow $\gamma^B_t$ on $T^*B$ with respect to the connection
$\pi_B^*\omega$.
\end{ex}

\begin{ex}
Now suppose that a fibration $f: M\to B$ as above is the orthonormal
frame bundle $F(B)\to B$ of the Riemannian manifold $B$. So, for any
$x\in B$, the fiber $F(B)_x$ consists of all orthonormal frames
$(v_1, v_2, \ldots, v_q)$ in $T_xB$. It is a principal bundle with
structure group $O(q)$. The Riemannian metric on $B$ gives rise to a
natural (Levi-Civita) connection on $f: F(B)\to B$. Fix a
bi-invariant Riemannian metric on $O(q)$ and consider the
corresponding Kaluza-Klein metric on $F(B)$.

By (\ref{e:nf}), it follows that
\[
N^*\cF\cong \{ ((v_1, v_2, \ldots, v_q),\xi)\in F(B)_x\times T_x^*B
: x\in B\}.
\]
For any $((v_1, v_2, \ldots, v_q),\xi)\in N^*\cF$, the action of the
transverse geodesic flow $\gamma^M_t$ is described as
\[
\gamma^M_t((v_1, v_2, \ldots, v_q),\xi) = ((v_1(t), v_2(t), \ldots,
v_q(t)),\xi(t)),
\]
where $\xi(t)=\gamma^B_t(\xi)$ and, for any $i=1,\ldots,q$, the
vector $v_i(t)$ is obtained by the parallel transport of $v_i$ along
the geodesic $\{\pi(\gamma^B_\tau(\xi)) :\tau\in [0,t]\}$ with
respect to the Levi-Civita connection on $TB$. Since
$\xi(t)=\gamma^B_t(\xi)$ can obtained by the parallel transport of
$\xi$ along the geodesic $\{\pi(\gamma^B_\tau(\xi)) :\tau\in
[0,t]\}$ with respect to the Levi-Civita connection on $T^*B$, the
transverse geodesic flow $\gamma^M_t$ has $q$ first integrals
$I_1,I_2,\ldots, I_q\in C^\infty(N^*\cF)$ given by
\[
I_j((v_1, v_2, \ldots, v_q),\xi)=\xi(v_j), \quad j=1,\ldots, q.
\]

There is a natural global right action of the group $SO(q)$ in the
fibers of the bundle
\[
N^*\cF\to T^*B, \quad ((v_1, v_2, \ldots, v_q),\xi)\in N^*\cF
\mapsto \xi\in T^*B.
\]
For every orthogonal matrix $A=(a_{ij})\in SO(q)$
and any $((v_1, v_2, \ldots, v_q),\xi)\in N^*\cF$ we put
\[
A ((v_1, v_2, \ldots, v_q),\xi)=\left((\sum_{i=1}^q v_ia_{i1},
\sum_{i=1}^q v_ia_{i2}, \ldots, \sum_{i=1}^q v_ia_{iq}),\xi\right).
\]
This action obviously commutes with the transverse geodesic flow
$\gamma^M_t$. Moreover, we have
\[
I_j(A(v_1, v_2, \ldots, v_q),\xi)=\sum_{i=1}^q I_i((v_1, v_2,
\ldots, v_q),\xi) a_{ij}.
\]
Therefore, the restrictions of the transverse geodesic flow
$\gamma^M_t$ to the level sets $(N^*\cF)_c$ defined by
\[
I_j((v_1, v_2, \ldots, v_q),\xi)=c_j, \quad j=1,2,\ldots,q,
\]
are isomorphic for different values of $c=(c_1,c_2,\ldots,c_q)\in
\RR^q$. It is easy to see that, for any $c\in\RR^q$, $(N^*\cF)_c$
can be identified with the frame bundle $F(B)$, and, for
$c=(1,0,\ldots,0)$, the restriction of $\gamma^M_t$ to $(N^*\cF)_c$
is precisely the frame flow on $F(B)$ (see \cite{Ja-Str06} and
references therein).
\end{ex}

\section{Egorov's theorem}\label{s:Egorov}
Let $D\in\Psi^1(M,E)$ be a formally self-adjoint, transversally
elliptic operator such that $D^2$ has the scalar principal symbol
and the holonomy invariant transverse principal symbol. By
\cite{noncom}, the operator $D$ is essentially self-adjoint with
initial domain $C^\infty(M,E)$. Define an unbounded linear
operator $\langle D\rangle$ in the space $L^2(M,E)$ as
\[
\langle D\rangle=(D^2+I)^{1/2}.
\]
By the spectral theorem, the operator $\langle D\rangle $ is
well-defined as a positive, self-adjoint operator in $L^2(M,E)$.
It can be shown that $H^1(M,E)$ is contained in the domain of
$\langle D\rangle$ in $L^2(M,E)$.

By the spectral theorem, the operator $\langle
D\rangle^s=(D^2+I)^{s/2}$ is a well-defined positive self-adjoint
operator in $\cH=L^2(M,E)$ for any $s\in\RR$, which is unbounded
if $s>0$. For any $s\geq 0$, denote by $\cH^s$ the domain of
$\langle D\rangle^s$, and, for $s<0$, $\cH^s=(\cH^{-s})^*$. Put
also $\cH^{\infty}=\bigcap_{s\geq 0}\cH^s, \quad
\cH^{-\infty}=(\cH^{\infty})^*$. It is clear that $H^s(M,E)
\subset \cH^s$ for any $s\geqslant 0$ and $\cH^s \subset H^s(M,E)$
for any $s<0$.  In particular, $C^\infty(M,E) \subset \cH^s$ for
any $s$.

We say that a bounded operator $A$ in $\cH^{\infty}$ belongs to
$\cL(\cH^{-\infty},\cH^{\infty})$ (resp.
$\cK(\cH^{-\infty},\cH^{\infty})$), if, for any $s$ and $r$, it
extends to a bounded (resp. compact) operator from $\cH^s$ to
$\cH^r$, or, equivalently, the operator $\langle
D\rangle^rA\langle D\rangle^{-s}$ extends to a bounded (resp.
compact) operator in $L^2(M,E)$. It is easy to see that
$\cL(\cH^{-\infty},\cH^{\infty})$ is a involutive subalgebra in
$\cL(\cH)$ and $\cK(\cH^{-\infty},\cH^{\infty})$ is its ideal. We
also introduce the class $\cL^1(\cH^{-\infty},\cH^{\infty})$,
which consists of all operators from
$\cK(\cH^{-\infty},\cH^{\infty})$ such that, for any $s$ and $r$,
the operator $\langle D\rangle^rA\langle D\rangle^{-s}$  is a
trace class operator in $L^2(M,E)$. It should be noted that any
operator $K$ with the smooth kernel belongs to
$\cL^1(\cH^{-\infty},\cH^{\infty})$.

By the spectral theorem, the operator $\langle D\rangle$ defines a
strongly continuous group $e^{it\langle D\rangle}$ of bounded
operators in $L^2(M,E)$. Consider a one-parameter group $\Phi_t$
of $\ast$-automorphisms of the algebra ${\mathcal L}(L^2(M,E))$
defined by
\begin{displaymath}
\Phi_t(T)=e^{i t\langle D\rangle}Te^{-i t\langle D\rangle}, \quad
T\in {\mathcal L}(L^2(M,E)), \quad t\in \RR.
\end{displaymath}

The main result of the paper is the following theorem.

\begin{thm}
\label{Egorov} Let $D\in\Psi^1(M,E)$ be a formally self-adjoint,
transversally elliptic operator such that $D^2$ has the scalar
principal symbol and the holonomy invariant transverse principal
symbol. For any $K\in \Psi^{m,-\infty}(M,{\mathcal F},E)$, there
exists an operator $K(t)\in\Psi^{m,-\infty}(M,{\mathcal F},E)$ such
that $\Phi_t(K)-K(t), t\in \RR,$ is a smooth family of operators of
class $\cL^1(\cH^{-\infty},\cH^{\infty})$.

Moreover, if $k\in S^m(G_{{\mathcal F}_N},{\mathcal
L}(\pi^*E)\otimes |T{\mathcal G}_N|^{1/2})$ is the principal symbol
of $K$, then the principal symbol $k_t\in S^m(G_{{\mathcal
F}_N},{\mathcal L}(\pi^*E)\otimes |T{\mathcal G}_N|^{1/2})$ of the
operator $K(t)$ is given by
\begin{equation}\label{e:ad}
k_t=\operatorname{Ad}(\alpha_t)^*(k),
\end{equation}
where $\operatorname{Ad}(\alpha_t)^*$ is the transverse
bicharacteristic flow of the operator $\langle D\rangle $.
\end{thm}

\begin{proof}
Let $\cL(\cD'(M,E),\cH^{\infty})$ (resp.
$\cL(\cH^{-\infty},C^{\infty}(M,E))$) be the space of all bounded
operators from $\cD'(M,E)$ to $\cH^{\infty}$ (resp. from
$\cH^{-\infty}$ to $C^{\infty}(M,E)$). Since any operator from
$\Psi^{-N}(M,E)$ with $N>\dim M$ is a trace class operator in
$L^2(M,E)$, one can easily see that
$\cL(\cD'(M,E),\cH^{\infty})\subset \cL^1(\cH^{-\infty},
\cH^{\infty})$ and $\cL(\cH^{-\infty},C^{\infty}(M,E))\subset
\cL^1(\cH^{-\infty}, \cH^{\infty})$.

As shown in \cite{mpag}, the operator $\langle
D\rangle=(D^2+I)^{1/2}$ can be written as $\langle D\rangle=P+R$,
where $P\in \Psi^1(M,E)$ is a self-adjoint, elliptic operator with
the positive, scalar principal symbol and the holonomy invariant
transversal principal symbol, and, for any $K\in
\Psi^{*,-\infty}(M,{\mathcal F},E)$,
$KR\in\cL(\cH^{-\infty},C^{\infty}(M,E))$ and $R K\in
\cL(\cD'(M,E),\cH^{\infty})$.

Denote by $e^{itP}$ the strongly continuous group of bounded
operators in $L^2(M,E)$ generated by the elliptic operator $i P$.
For $K\in \Psi^{m,-\infty}(M,{\mathcal F},E)$, let
$\Phi_t^{P}(K)=e^{itP}Ke^{-itP}$. It is shown in \cite{mpag} that
the operator $\Phi_t^{P}(K)=e^{itP}Ke^{-itP}$ is in
$\Psi^{m,-\infty}(M,{\mathcal F},E)$, and $\Phi_t(K)-\Phi_t^{P}(K),
t\in \RR,$ is a smooth family of operators of class
$\cL^1(\cH^{-\infty},\cH^{\infty})$. So we can take
$K(t)=\Phi_t^{P}(K)$. It remains to compute the principal symbol of
$\Phi_t^{P}(K)$.

Without loss of generality, one can assume that the elliptic
extension $\tilde p$ of $p$ introduced in Section~\ref{s:flow} to
define the transverse bicharacteristic flow coincides with the
principal symbol of $P$. We have
\[
\frac{d}{d t}\Phi^P_t(K)=[iP,\Phi^P_t(K)], \quad t\in \RR, \quad
\Phi^P_0(K)=K.
\]
Recall (cf. (\ref{e:subflow})) that the function $k_t$ given by
(\ref{e:ad}) satisfies the following equation
\begin{equation}\label{e:kt}
\frac{d}{d t}k_t=\nabla_{\cH_{p}}k_t.
\end{equation}
Let $K_0(t)$ be any operator from $\Psi^{m,-\infty}(M,{\mathcal
F},E)$ with the principal symbol $k_t$. Then, by (\ref{e:comm2}) and
(\ref{e:kt}), it follows that
\[
\frac{d}{d t}K_0(t)=[iP,K_0(t)]+R(t), \quad t\in \RR, \] \[
K_0(0)=K+R_0.
\]
where $R(t)\in \Psi^{m-1,-\infty}(M,{\mathcal F},E), t\in \RR,$ and
$R_0\in \Psi^{m-1,-\infty} (M,{\mathcal F},E)$. It is easy to see
that
\[
K_0(t)-\Phi_t^{P}(K)=\int_0^t\Phi_{t-\tau}(R(\tau))d\tau+\Phi_t(R_0),
\]
that immediately implies that $K_0(t)-\Phi_t^{P}(K)\in
\Psi^{m-1,-\infty} (M,{\mathcal F},E)$.
\end{proof}

\section{Preliminaries on transverse Dirac operators} Let $M$ be a compact manifold
equipped with a Riemannian foliation $\cF$ of even codimension $q$
and $\cE$ a Hermitian vector bundle over $M$ equipped with a
leafwise flat unitary connection $\nabla^\cE$. Suppose that $g_M$ is
a bundle-like metric on $M$.

As above, let $T^H_xM=T_x{\cF}^{\bot}$. Let $P_H$ (resp. $P_F$)
denotes the orthogonal projection operator of $TM=T^HM\oplus T\cF$
on $T^HM$ (resp. $T\cF$). There is the canonical flat connection
$\stackrel{\circ}{\nabla}$ in $T^HM$, defined along the leaves of
$\cF$ (the Bott connection) given by
\[
{\stackrel{\circ}\nabla}_X N= P_H[X,N],\quad X\in
C^\infty(M,T\cF),\quad N\in C^\infty(M,T^HM).
\]
Denote by $\nabla^L$ the Levi-Civita connection defined by $g_M$.
The following formulas define a connection $\nabla$ in $T^HM$
(called the transverse Levi-Civita connection):
\begin{equation}\label{e:adapt}
\begin{aligned}
\nabla_XN &=P_H[X,N],\quad X\in C^\infty(M,T\cF),\quad N\in
C^\infty(M,T^HM)\\ \nabla_XN&=P_H\nabla^L_XN,\quad X\in
C^\infty(M,T^HM),\quad N\in C^\infty(M,T^HM).
\end{aligned}
\end{equation}
It turns out that $\nabla$ depends only on the transverse part of
the metric $g_M$ and preserves the inner product of $T^HM$. This
connection will be called the transverse Levi-Civita connection.

Denote by $\cR$ the integrability tensor (or curvature) of $T^HM$.
It is the $2$-form on $T^HM$ with values in $T\cF$ given by
\[
\cR_x(f_1,f_2)=-P_F[\tilde{f}_1,\tilde{f}_2](x), \quad f_1, f_2\in
T^H_xM,
\]
where, for any $f\in T^H_xM$, $\tilde{f}\in C^\infty(M,T^HM)$
denotes any infinitesimal transformation of $\cF$, which coincides
with $f$ at $x$.

Since the Levi-Civita connection $\nabla^L$ has no torsion, for any
$f_1, f_2\in C^\infty(M,T^HM)$, we have
\begin{equation}
\label{e:R} \nabla_{f_1}f_2-\nabla_{f_2}f_1=P_H([f_1,f_2])=
[f_1,f_2]+\cR(f_1,f_2).
\end{equation}

Let $\omega_\cF$ denote the leafwise Riemannian volume form of
$\cF$. Let $f\in T^H_xM$ and let $\tilde{f}\in C^\infty(M,T^HM)$
denote any infinitesimal transformation of $\cF$, which coincides
with $f$ at $x$. The local flow generated by $\tilde{f}$ preserves
the foliation and gives rise to a well-defined action on
$\Lambda^pT^*\cF$. The mean curvature vector field $\tau\in
C^\infty(M,T^HM)$ of $\cF$ is defined by the identity
\[
L_{\tilde{f}}\omega_\cF=g_M(\tau, \tilde{f}) \omega_\cF
\]
If $e_1,e_2,\ldots,e_p$ is a local orthonormal frame in $T\cF$, then
\[
\tau=\sum_{i=1}^pP_H(\nabla^L_{e_i}e_i).
\]

Assume that $\cF$ is transversely oriented and the normal bundle $Q$
is spin. Thus the $SO(q)$ bundle $O(Q)$ of oriented orthonormal
frames in $Q$ can be lifted to a $Spin(q)$ bundle $O'(Q)$ so that
the projection $O'(Q)\to O(Q)$ induces the covering projection
$Spin(q)\to SO(q)$ on each fiber.

Let $F(Q), F_+(Q), F_-(Q)$ be the bundles of spinors
\[
F(Q)=O'(Q)\times_{Spin(q)}S, \quad
F_\pm(Q)=O'(Q)\times_{Spin(q)}S_\pm .
\]
Denote by $Cl(Q_x)$ the Clifford algebra of $Q_x$, $x\in M$. Recall
that, relative to an orthonormal basis $\{f_1,f_2,\ldots,f_q\}$ of
$Q_x$, $Cl(Q_x)$ is the complex algebra generated by $1$ and
$f_1,f_2,\ldots,f_q$ with relations
\[
f_\alpha f_\beta+f_\beta f_\alpha=-2\delta_{\alpha\beta}, \quad
\alpha, \beta=1,2,\ldots,q.
\]
Since $\dim Q=q$ is even $\End F(Q)$ is as a bundle of algebras
over $M$ isomorphic to the Clifford bundle $Cl(Q)$. The action of
an element $a\in Cl(Q)$ on $F(Q)$ will be denoted by $c(a)$.

The transverse Levi-Civita connection $\nabla$ lifts to a connection
$\nabla^{F(Q)}$ on the holonomy equivariant vector bundle $F(Q)$,
whose restriction to $T\cF$ coincides with the Bott connection. It
can be easily seen that $\nabla^{F(Q)}$ is a Clifford connection,
that is, for any $f\in T^HM$ and $X\in T^HM$, we have
\[
[\nabla^{F(Q)}_f, c(X)]=c(\nabla_fX).
\]
Let
\[
\nabla^{F(Q)\otimes \cE}=\nabla^{F(Q)}\otimes 1 + 1\otimes
\nabla^{\cE}
\]
be the corresponding connection on $F(Q)\otimes \cE$.

We will identify the bundle $Q$ and $Q^*$ by means of the metric
$g_M$ and define the operator $D^\prime_\cE$ acting on the sections
of $F(Q)\otimes \cE$ as the composition
\begin{multline*}
C^\infty(M,F(Q)\otimes \cE)\stackrel{\nabla^{F(Q)\otimes
\cE}}{\longrightarrow} C^\infty(M,Q^*\otimes F(Q)\otimes \cE) \\
= C^\infty(M,Q\otimes F(Q)\otimes \cE) \stackrel{c\otimes
1}{\longrightarrow} C^\infty(M,F(Q)\otimes \cE).
\end{multline*}
This operator is odd with respect to the $\ZZ_2$-grading
$F(Q)\otimes \cE=(F_+(Q)\otimes \cE)\oplus (F_-(Q)\otimes \cE)$. If
$f_1,\ldots,f_q$ is a local orthonormal frame for $T^HM$, then
\[
D^\prime_\cE=\sum_{\alpha=1}^q(c(f_\alpha)\otimes
1)\nabla^{F(Q)\otimes \cE}_{f_\alpha}.
\]

Denote by $(\cdot,\cdot)_x$ the inner product in the fiber
$(F(Q)\otimes \cE)_x$ over $x\in M$. Then the inner product in
$L^2(M,F(Q)\otimes \cE)$ is given by the formula
\[
(s_1, s_2)=\int_M (s_1(x),s_2(x))_x\omega_M, \quad s_1, s_2\in
L^2(M,F(Q)\otimes \cE),
\]
where $\omega_M=\sqrt{\det g}\,dx$ denotes the Riemannian volume
form on $M$. In the following lemma, we compute the formal adjoint
$(D^\prime_\cE)^*$ of $D^\prime_\cE$ (see also \cite{G-K91} and
references therein).

\begin{lem}
We have
\[
(D^\prime_\cE)^*=D^\prime_\cE-c(\tau).
\]
\end{lem}

\begin{proof}
For any $s_1, s_2\in C^\infty(M,F(Q)\otimes \cE)$, we have
\begin{multline*}
(D^\prime_\cE s_1, s_2) \\ = \sum_{\alpha=1}^q \left( - \int_M
f_\alpha [(s_1, (c(f_\alpha)\otimes 1) s_2)_x]\omega_M +(s_1,
\nabla^{F(Q)\otimes \cE}_{f_\alpha} (c(f_\alpha)\otimes 1) s_2)
\right) \\ = \sum_{\alpha=1}^q \left( - \int_M f_\alpha [(s_1,
(c(f_\alpha)\otimes 1) s_2)_x]\omega_M +(s_1, (c(\nabla_{f_\alpha}
f_\alpha)\otimes 1) s_2)\right)\\ +(s_1, D^\prime_\cE s_2).
\end{multline*}
Recall that, by the divergence theorem, for any vector field $X$ on
$M$ and $a\in C^\infty(M)$, we have
\[
\int_M X(a)(x)\omega_M=-\int_M \operatorname{div}(X) \cdot
a(x)\omega_M.
\]
Let $e_1,e_2,\ldots,e_p$ be a local orthonormal frame in $T\cF$.
Then the divergence $\operatorname{div}(X)$ of $X$ is given by the
formula
\begin{equation}\label{e:div}
\operatorname{div} (X)= \sum_{k=1}^p g_M(e_k,
\nabla_{e_k}X)+\sum_{\beta=1}^q g_M(f_\beta, \nabla_{f_\beta}X).
\end{equation}
In particular, it is easy to see that
\[
\operatorname{div} (f_\alpha)= - g_M(\tau+
\sum_{\beta=1}^q\nabla_{f_\beta}f_\beta, f_\alpha).
\]
Using the divergence theorem, we easily get
\begin{multline*}
\sum_{\alpha=1}^q \left( - \int_M f_\alpha [(s_1,
(c(f_\alpha)\otimes 1)s_2)_x]\omega_M +(s_1, (c(\nabla_{f_\alpha} f_\alpha)\otimes 1) s_2)\right)\\
= - (s_1, (c(\tau)\otimes 1) s_2),
\end{multline*}
that completes the proof.
\end{proof}

By this lemma, the operator
\[
D_\cE=D^\prime_\cE-\frac12 c(\tau)= \sum_{\alpha=1}^q
(c(f_\alpha)\otimes 1)\left(\nabla^{F(Q)\otimes
\cE}_{f_\alpha}-\frac12 g_M(\tau, f_\alpha) \right)
\]
is self-adjoint. This operator will be called the transverse Dirac
operator. It was introduced in \cite{G-K91} (see also
\cite{G-K91a,G-K93} and references therein).

We will use the Riemannian volume form $\omega_M$ to identify the
half-densities bundle with the trivial one. So the action of $D_\cE$
on half-densities is defined by
\[
D_\cE (u|\omega_M|^{1/2})=(D_\cE u)|\omega_M|^{1/2}, \quad u\in
C^\infty(M,F(Q)\otimes \cE).
\]

\section{The transverse signature operator}
In this section, we will discuss a particular example of a
transverse Dirac operator given by the transverse signature
operator.

As above, let $(M,{\mathcal F})$ be a compact Riemannian foliated
manifold equipped with a bundle-like metric $g_M$.

\begin{lem}\label{l:dH}
Let $f_1, f_2, \ldots, f_q$ be a local orthonormal basis of $T^HM$
and $f^*_1, f^*_2, \ldots, f^*_q$ be the dual basis of $T^HM^*$.
Then on $C^\infty(M,\Lambda T^HM^{*})$ we have
\begin{gather*}
d_H=\sum_{\alpha=1}^q\varepsilon_{f^*_\alpha} \nabla_{f_\alpha},\\
d^*_H=-\sum_{\alpha=1}^q i_{f_\alpha} \nabla_{f_\alpha} + i_\tau.
\end{gather*}
\end{lem}

\begin{proof}
Denote $d^\prime_H=\sum_{\alpha=1}^q\varepsilon_{f^*_\alpha}
\nabla_{f_\alpha}$. Then the operators $d_H$ and $d^\prime_H$
satisfy the Leibniz rule and, clearly, coincide on functions. It
remains to show that they agree on the space $C^\infty(M,T^HM^*)$ of
transverse one-forms. Using an explicit formula for $d_H$ and
(\ref{e:R}), for any $\omega \in C^\infty(M,T^HM^*)$ and for any
$U,V\in C^\infty(M,T^HM)$, we get
\[
\begin{split}
d_H\omega(U,V) & =U[\omega(V)]-V[\omega(U)]-\omega(P_H[U,V])\\ & =
\nabla_U\omega(V)+\omega(\nabla_UV)-\nabla_V\omega(U)-\omega(\nabla_VU)-\omega(P_H[U,V])\\
& = \nabla_U\omega(V)-\nabla_V\omega(U).
\end{split}
\]

Now, since $U=\sum_\alpha \langle f^*_\alpha, U\rangle f_\alpha$
and $V=\sum_\alpha \langle f^*_\alpha, V\rangle f_\alpha$, we
obtain
\[
\nabla_U\omega(V)-\nabla_V\omega(U)=\sum_{\alpha}(\langle
f^*_\alpha, U\rangle \nabla_{f_\alpha}\omega(V)-\langle
f^*_\alpha, V\rangle \nabla_{f_\alpha}\omega(U))
=d^\prime_H\omega(U,V),
\]
that proves the first equality.

The second equality can be easily derived from the first one, if we
take the adjoints and use the divergence theorem.
\end{proof}

To represent the transverse signature operator $D_H=d_H+d_H^*$ as a
transverse Dirac operator, we take $\cE=F(Q)^*$. By
Lemma~\ref{l:dH}, we have the following formula for the
corresponding transverse Dirac operator $D_{F(Q)^*}$:
\[
\begin{split}
D_{F(Q)^*}&
=\sum_{\alpha=1}^q(\varepsilon_{f^*_\alpha}-i_{f_\alpha})
\nabla_{f_\alpha}-\frac12(\varepsilon_{\tau^*}-i_{\tau})\\ & =
d_H+d^*_H-\frac12(\varepsilon_{\tau^*}+i_{\tau}).
\end{split}
\]
So we see that the transverse signature operator $D_H$ coincides
with the transverse Dirac operator $D_{F(Q)^*}$ if and only if
$\tau=0$, that is, all the leaves are minimal submanifolds.

\begin{ex}
Consider a foliation $\cF$ on a compact manifold $M$ given by the
fibers of a principal $K$-bundle $f: M\to B$ with connection, where
$K$ is a compact group. Fix a Riemannian metric on $B$ and a
bi-invariant metric on $K$, and consider the corresponding
Kaluza-Klein metric on $M$.

For any irreducible unitary representation $\rho$ of $K$ in a vector
space $W_\rho$, consider the associated Hermitian vector bundle
$E_\rho=M\times_\rho W_\rho$ over $B$. It is well-known that there
is a natural identification of the space $C^\infty(B, E_\rho)$ with
the space $F_\rho$ of smooth functions $f : M\to W_\rho$ satisfying
$f(x\cdot k)=\rho(k)^{-1}f(x)$ for any $x\in M$ and $k\in K$. We
denote by $C^\infty(M)_\rho$ the isotypical component of $\rho$ in
$C^\infty(M)$. So we have
\[
C^\infty(M)=\bigoplus_{\rho\in \hat{K}} C^\infty(M)_\rho.
\]
The following lemma is a generalization of the usual Peter-Weyl
theorem to bundles (see, for instance, \cite[Lemma
5.3]{Gui-Uribe86}).

\begin{lem}
The mapping
\[
J_\rho : C^\infty(B, E_\rho)\otimes W^*_\rho\to C^\infty(M), \quad
f\otimes \eta \mapsto f_\eta,
\]
where
\[
f_\eta(x)=\sqrt{\frac{\dim W_\rho}{\operatorname{vol} K}}\eta(f(x)),
\quad x\in M,
\]
is a unitary isomorphism onto $C^\infty(M)_\rho$, which is
$K$-equivariant with respect to the representation $1\otimes \rho^*$
on $C^\infty(B, E_\rho)\otimes W^*_\rho$.
\end{lem}

Next, the transverse de Rham differential $d_H$ commutes with the
natural action of $K$ on $C^{\infty}(M, \Lambda^*T^HM^{*})$.
Therefore, $d_H$ maps $C^{\infty}(M, \Lambda^*T^HM^{*})_\rho$ to
$C^{\infty}(M, \Lambda^*T^HM^{*})_\rho$. Let $ \nabla^{E_\rho} :
C^\infty(B, \Lambda^*T^*B\otimes E_\rho)\to C^\infty(B,
\Lambda^*T^*B\otimes E_\rho)$ is the exterior covariant derivative
associated with the connection. By definition (see, for instance
\cite{K-N}), under the isomorphism $J_\rho$, the restriction of
$d_H$ to $C^{\infty}(M, \Lambda^*T^HM^{*})_\rho$ corresponds to the
operator $\nabla^{E_\rho}\otimes I_{W_\rho^*}$ on $C^\infty(B,
\Lambda^*T^*B \otimes E_\rho)\otimes W^*_\rho$. Since the
isomorphism $J_\rho$ is unitary, the similar statement holds for
$d^*_H$.

Thus, we have the commutative diagram
\[
  \begin{CD}
C^\infty(B, \Lambda^*T^*B\otimes E_\rho)\otimes W^*_\rho @>
D^{E_\rho} \otimes
I_{W_\rho^*}>> C^\infty(B, \Lambda^*T^*B\otimes E_\rho)\otimes W^*_\rho\\
@VVV @VVV
\\ C^{\infty}(M, \Lambda^*T^HM^{*})_\rho @>D_H>> C^{\infty}(M, \Lambda^*T^HM^{*})_\rho
  \end{CD}
\]
where $D^{E_\rho}=\nabla^{E_\rho}+(\nabla^{E_\rho})^*$ is the
twisted signature operator on $B$ with coefficients in the vector
bundle $E_\rho$. It shows that, in this case, the transverse
signature operator $D_H=d_H+d_H^*$ decomposes into a direct sum of
twisted signature operators on the base $B$ with coefficients in
vector bundles associated with irreducible representations of $K$.
\end{ex}

\section{The subprincipal symbol of a transverse Dirac operator}
In this section we compute the transverse bicharacteristic flow of
transverse Dirac operators. For this, we will use the following
fact (see, for instance, \cite[Proposition 4.3.1]{Duistermaat}).

\begin{thm}\label{t:sub}
Let $P\in \Psi^m(X)$ be a properly supported pseudodifferential
operator on a smooth manifold $X$. For any $a\in C^\infty(X,
|TX|^{1/2})$ and for any real-valued function $\phi\in C^\infty(X)$
we have
\begin{multline*}
e^{-is\phi(x)}P(e^{is\phi}a)(x)=s^mp_m(x,d\phi(x))\cdot a(x)\\ +
s^{m-1}\left(p_{sub}(x,d\phi(x))\cdot
a(x)+\frac{1}{i}(\cL_va)(x)\right)+O(s^{m-2}), \quad s\to \infty,
\end{multline*}
where $v$ is a vector field on $X$:
\[
v(x)=\sum_j \frac{\partial p_m}{\partial \xi_j} (x,
d\phi(x))\frac{\partial }{\partial x_j}=\pi_*(X_p(x, d\phi(x))),
\]
$X_p$ is the Hamiltonian vector field of $p_m$ on $T^*X$,
$\pi_*(X_p(x,\xi))\in T_xX$ is the image of $X_p(x,\xi)\in
T_{(x,\xi)}(T^*X)$ under the projection $\pi:T^*X\to X$.
\end{thm}

Here $\cL_v$ denotes the Lie derivative along $v$, acting on
half-densities: for any $f\in C^\infty(X)$, we have
\[
\cL_v(f|\omega_X|^{1/2})=v(f)|\omega_X|^{1/2}+\frac{1}{2}\operatorname{div}
v\cdot f|\omega_X|^{1/2}.
\]
This theorem remains to be true for operators acting on sections of
a vector bundle $E$ over $X$ locally, that is, if we fix a
trivialization of $E$ over some open subset of $X$.

Let $M$ be a compact manifold equipped with a Riemannian foliation
$\cF$ of even codimension $q$, $\cE$ a Hermitian vector bundle over
$M$ equipped with a leafwise flat unitary connection $\nabla^\cE$,
$g_M$ a bundle-like metric on $M$ and $D_\cE$ the associated
transverse Dirac operator.

\begin{thm}\label{t:trsym}
The principal symbol of $D^2_\cE$ is given by
\begin{equation}\label{e:principal}
a_2(x,\xi)=g^M(P^H(x,\xi), P^H(x,\xi)),\quad (x,\xi)\in T^*M,
\end{equation}
where $g^M$ is the induced metric on $T^*M$, $P^H: T^*M\to T^HM^*$
is the orthogonal projection.
\end{thm}

\begin{proof}
Let $f_1,\ldots, f_q$ be a local orthonormal basis of $T^HM$, which
consists of infinitesimal transformations of $\cF$. For any $a\in
C^\infty(M,F(Q)\otimes \cE)$ and for any real-valued function
$\phi\in C^\infty(M)$ we have
\begin{multline}\label{e:tau}
e^{-is\phi(x)}D^2_\cE(e^{is\phi}a)(x) \\ = \left(\sum_{\alpha=1}^q
(c(f_\alpha)\otimes 1)(\nabla^{F(Q)\otimes
\cE}_{f_\alpha}-\frac12g_M(\tau, f_\alpha) +is
f_\alpha(\phi))\right)^2.
\end{multline}

The terms of order $s^2$ in (\ref{e:tau}) are
\[
s^2\sum_{\alpha=1}^q (f_\alpha(\phi))^2=s^2\sum_{\alpha=1}^q \langle
d\phi, f_\alpha \rangle^2.
\]
Therefore, the principal symbol of $D_\cE^2$ is
\[
a_2(x,\xi)=\sum_{\alpha=1}^q \langle \xi, f_\alpha(x) \rangle^2,
\quad (x,\xi)\in T^*M,
\]
that completes the proof.
\end{proof}

It is easy to see that the 1-jets of the functions $\sqrt{a_2}$,
where $a_2$ is the principal symbol of $D^2_\cE$, and
\[
p(x,\xi)=\sqrt{g^M((x,\xi), (x,\xi))},\quad (x,\xi)\in T^*M,
\]
coincide on $N^*\cF$. So we can further work with the elliptic
symbol $p$.

The Hermitian connection $\nabla^{F(Q)\otimes \cE}$ determines
uniquely a Hermitian partial connection $\tilde{\nabla}_{X_p}$ along
the Hamiltonian vector field $X_p$ on $\pi^*(F(Q)\otimes \cE)$,
which satisfies
\[
\left(\tilde{\nabla}_{X_p}\pi^*s\right)(\nu)=\nabla^{F(Q)\otimes
\cE}_{\pi_*(X_p(\nu))}s(\pi(\nu)), \quad s\in C^\infty(M,F(Q)\otimes
\cE),
\]
where $\pi^*s\in C^\infty(\tilde{N}^*\cF, \pi^*(F(Q)\otimes \cE))$
is the pull back of a section $s\in C^\infty(M,F(Q)\otimes \cE)$
under the projection $\pi : \tilde{N}^*\cF \to M$.

If we fix a local orthonormal basis $f_1,\ldots, f_q$ of $T^HM$,
which consists of infinitesimal transformations of $\cF$, and a
local trivialization of $F(Q)\otimes \cE$ and write
$\nabla^{F(Q)\otimes \cE}_{f_\alpha}=f_\alpha+B(f_\alpha)$ with some
matrix-valued one form $B$, then, for $s\in C^\infty(\tilde{N}^*\cF,
\pi^*(F(Q)\otimes \cE))$, we have
\[
(\tilde{\nabla}_{X_p}s)(\nu)=
X_{p}s(\nu)+\|\nu\|^{-1}\sum_{\alpha=1}^q \langle \nu,f_\alpha
\rangle B(f_\alpha)s(\nu), \quad \nu\in N^*\cF.
\]

The geometric meaning of this partial connection is as follows.
Recall that $X_p$ generates the geodesic flow $f_t$ on
$\tilde{N}^*\cF$. For any $\nu \in \tilde{N}^*\cF$, the projection
of the orbit $\cO_\nu=\{f_t(\nu), t\in \RR\}$ to $M$ is the geodesic
$\gamma_\nu$, passing through $x=\pi(\nu)$. Then the parallel
transport of $v\in \pi^*(F(Q)\otimes \cE)_\nu$ along $\cO_\nu$ with
respect to the connection $\tilde{\nabla}_{X_p}$ coincides with the
parallel transport of $v$ considered as an element of $(F(Q)\otimes
\cE)_x$ along the geodesic $\gamma_\nu$ with respect to the
connection $\nabla^{F(Q)\otimes \cE}$.

\begin{thm}
The subprincipal symbol of $\langle D_\cE\rangle$ considered as a
partial connection $\nabla_{sub}(\langle D_\cE\rangle)$  on
$\pi^*(F(Q)\otimes \cE)$ coincides with $\tilde{\nabla}_{X_p}$.
\end{thm}

\begin{proof}
The terms of order $s$ in (\ref{e:tau}) are
\begin{multline*}
i[(\sum_{\alpha=1}^q(c(f_\alpha)\otimes
1){f_\alpha}(\phi))(\sum_{\beta=1}^q(c(f_\beta)\otimes
1)(\nabla^{F(Q)\otimes \cE}_{f_\beta}-\frac12g_M(\tau,
f_\beta)))\\ +(\sum_{\beta=1}^q(c(f_\beta)\otimes
1)(\nabla^{F(Q)\otimes \cE}_{f_\beta}-\frac12g_M(\tau, f_\beta)))
(\sum_{\alpha=1}^q(c(f_\alpha)\otimes 1){f_\alpha}(\phi))]\\ = i
\sum_{\alpha,\beta} ((c(f_\alpha)c(f_\beta) +
c(f_\beta)c(f_\alpha))\otimes 1) {f_\alpha}(\phi)
(\nabla^{F(Q)\otimes \cE}_{f_\beta}-\frac12g_M(\tau, f_\beta)) \\
+ i \sum_{\alpha,\beta}(c(f_\beta)c(\nabla_{f_\beta}
f_\alpha)\otimes 1){f_\alpha}(\phi)+ i \sum_{\alpha,\beta}
(c(f_\beta)c(f_\alpha)\otimes 1) f_\beta f_\alpha(\phi)\\
=I_1+I_2+I_3.
\end{multline*}
For the first term, we easily get
\[
I_1 =-2i \sum_{\alpha} {f_\alpha}(\phi)\nabla^{F(Q)\otimes
\cE}_{f_\alpha}-i\tau(\phi).
\]
Let $\nabla_{f_\alpha}f_\beta=\sum_\gamma
a^\gamma_{\alpha\beta}f_\gamma$. Since $\nabla$ is compatible with
the metric, we have $a^\gamma_{\alpha\beta}
=-a^\beta_{\alpha\gamma}$. Thus we get
\begin{multline*}
\begin{aligned}
I_2 = &
\frac{i}{2}\sum_{\alpha,\beta}[(c(f_\alpha)c(\nabla_{f_\alpha}
f_\beta)\otimes 1){f_\beta}(\phi)+(c(f_\beta)c(\nabla_{f_\beta}
f_\alpha)\otimes 1){f_\alpha}(\phi)]\\ =&
\frac{i}{2}\sum_{\alpha,\beta,\gamma} [a^\gamma_{\alpha\beta}
(c(f_\alpha)c(f_\gamma)\otimes 1){f_\beta}(\phi)
+a^\gamma_{\alpha\beta} (c(f_\beta)c(f_\gamma)\otimes
1){f_\alpha}(\phi)]\\ =& - \frac{1}{2}\sum_{\alpha,\beta,\gamma}
[a^\beta_{\alpha\gamma} (c(f_\alpha)c(f_\gamma)\otimes
1){f_\beta}(\phi)+a^\alpha_{\beta\gamma}
(c(f_\beta)c(f_\gamma)\otimes 1){f_\alpha}(\phi)]
\\ =& - \frac{i}{2}\sum_{\alpha,\gamma} [
(c(f_\alpha)c(f_\gamma)\otimes 1){\nabla_{f_\alpha}
f_\gamma}(\phi) + \sum_{\beta,\gamma}
(c(f_\beta)c(f_\gamma)\otimes 1){\nabla_{f_\beta} f_\gamma}(\phi)]
\\ =& - i \sum_{\alpha,\beta} (c(f_\alpha)c(f_\beta)\otimes
1){\nabla_{f_\alpha} f_\beta}(\phi).
\end{aligned}
\end{multline*}
Finally, we have
\[
\begin{split}
I_3=& \frac{i}{2} \sum_{\alpha,\beta} ((c(f_\beta)c(f_\alpha)\otimes
1) f_\beta f_\alpha(\phi)+ (c(f_\alpha)c(f_\beta)\otimes 1)
f_\alpha f_\beta(\phi)) \\
 = & \frac{i}{2} \sum_{\alpha,\beta}
(c(f_\beta)c(f_\alpha)+c(f_\alpha)c(f_\beta))\otimes 1) f_\alpha
f_\beta(\phi)\\ & + \frac{i}{2} \sum_{\alpha,\beta}
(c(f_\beta)c(f_\alpha)\otimes 1) [f_\beta, f_\alpha](\phi)\\
 = & -i \sum_{\alpha} f^2_\alpha
(\phi)+ \frac{i}{2} \sum_{\alpha,\beta}
(c(f_\beta)c(f_\alpha)\otimes 1) (\nabla_{f_\beta} f_\alpha -
\nabla_{f_\alpha} f_\beta-\cR(f_\beta, f_\alpha)) (\phi),
\end{split}
\]
where we used the equality (\ref{e:R}).

From the last three identities, the terms of order $s$ are
\begin{multline*}
-2i \sum_{\alpha} {f_\alpha}(\phi)\nabla^{F(Q)\otimes
\cE}_{f_\alpha} -i\tau(\phi) -i \sum_{\alpha} f^2_\alpha (\phi)\\ +
i \sum_{\alpha} \nabla_{f_\alpha} f_\alpha (\phi) - \frac{i}{2}
\sum_{\alpha,\beta} (c(f_\beta)c(f_\alpha)\otimes 1)\cR(f_\beta,
f_\alpha)(\phi).
\end{multline*}
By Theorem~\ref{t:sub}, we have
\begin{multline*}
p_{sub}(x,d\phi(x))\cdot
a(x)|\omega_M|^{1/2}+\frac{1}{i}\cL_v(a|\omega_M|^{1/2})(x)\\ =(- i
\sum_{\alpha=1}^q f^2_\alpha(\phi)a - 2 i \sum_{\alpha=1}^q
f_\alpha(\phi) \nabla^{F(Q)\otimes \cE}_{f_\alpha}a -i\tau(\phi)a +
i
\sum_{\alpha=1}^q \nabla_{f_\alpha} f_\alpha(\phi)a \\
-\frac{1}{2} i \sum_{\alpha=1}^q \sum_{\beta=1}^q
c(f_\alpha)c(f_\beta) \cR(f_\alpha,f_\beta)(\phi)a)|\omega_M|^{1/2}.
\end{multline*}
Now compute the vector field $v$:
\[
\pi_*(X_p(x,\xi))=2\sum_{\alpha=1}^q \langle \xi, f_\alpha \rangle
f_\alpha, \quad \xi\in T^*M
\]
and
\[
v=2\sum_{\alpha=1}^q \langle d\phi(x), f_\alpha \rangle
f_\alpha=2\sum_{\alpha=1}^q f_\alpha(\phi) f_\alpha.
\]
Therefore, we have
\[
\cL_v(a|\omega_M|^{1/2})=2\sum_{\alpha=1}^q f_\alpha(\phi)
f_\alpha(a)|\omega_M|^{1/2}+\sum_{\alpha=1}^q
\operatorname{div}(f_\alpha(\phi) f_\alpha) \cdot a|\omega_M|^{1/2}.
\]
Let $e_1,e_2,\ldots,e_p$ be a local orthonormal frame in $T\cF$.
Using (\ref{e:div}), we easily compute
\[
\sum_{\alpha=1}^q \operatorname{div}(f_\alpha(\phi) f_\alpha)=
\sum_{\alpha=1}^q f^2_\alpha(\phi) - \tau(\phi) - \sum_{\beta=1}^q
\nabla_{f_\beta} f_\beta(\phi).
\]
Finally, if we write $\nabla^{F(Q)\otimes
\cE}_{f_\alpha}=f_\alpha+B(f_\alpha)$ with some matrix-valued one
form $B$, we get
\[
p_{sub}(x,d\phi(x))=- 2 i \sum_{\alpha=1}^q
f_\alpha(\phi)B(f_\alpha)
\\ -\frac{i}{2} \sum_{\alpha=1}^q \sum_{\beta=1}^q
c(f_\alpha)c(f_\beta) \cR(f_\alpha,f_\beta)(\phi),
\]
or, equivalently,
\[
p_{sub}(x,\xi)=- 2 i \sum_{\alpha=1}^q \langle \xi,f_\alpha \rangle
B(f_\alpha)\\ -\frac{i}{2} \sum_{\alpha=1}^q \sum_{\beta=1}^q
c(f_\alpha)c(f_\beta) \langle \xi, \cR(f_\alpha,f_\beta)\rangle.
\]

The transverse subprincipal symbol of $D_\cE^2$ is
\[
\sigma_{sub}(D_\cE^2)(\nu)=- 2 i \sum_{\alpha=1}^q \langle
\nu,f_\alpha \rangle B(f_\alpha), \quad \nu\in N^*\cF.
\]
Since the principal symbol of $D_\cE^2$ is scalar, the formula
proved in \cite{DG}
\[
\sigma_{sub}(\langle D_\cE\rangle)
=\frac{1}{2}\sigma(D_\cE^2)^{-\frac{1}{2}}\sigma_{sub}(D_\cE^2),
\]
continues to hold and the transverse subprincipal symbol of
$\langle D_\cE\rangle$ is
\[
\sigma_{sub}(\langle D_\cE\rangle)(\nu)=- i
\|\nu\|^{-1}\sum_{\alpha=1}^q \langle \nu,f_\alpha \rangle
B(f_\alpha), \quad \nu\in N^*\cF.
\]
Thus, the subprincipal symbol of $\langle D_\cE\rangle$ is a partial
connection $\nabla_{sub}(\langle D_\cE\rangle)$ on $\pi^*E$ along
the Hamiltonian vector field $X_p$ given by
\[
\nabla_{sub}(\langle
D_\cE\rangle)=X_{p}(\nu)+\|\nu\|^{-1}\sum_{\alpha=1}^q \langle
\nu,f_\alpha \rangle B(f_\alpha)=\tilde{\nabla}_{X_p(\nu)}, \quad
\nu\in N^*\cF,
\]
as desired.
\end{proof}

\section{The noncommutative geodesic flow}
\label{ncg}

As stated in \cite{noncom} (see also \cite{mpag}), any operator $D$,
satisfying the assumptions of Section~\ref{s:Egorov}, defines a
spectral triple in the sense of Connes' noncommutative geometry
\cite{Sp-view,Co-M}. In this setting, Theorem~\ref{Egorov} has a
natural interpretation in terms of the corresponding noncommutative
geodesic flow.

More precisely, consider spectral triples $({\mathcal A},{\mathcal
H},D)$ associated with a compact foliated Riemannian manifold
$(M,{\mathcal F})$ (see \cite{mpag} for more details):
\begin{enumerate}
\item The involutive algebra ${\mathcal A}$ is the algebra $\cinf_c(G,|T\cG|^{1/2})$;
\item The Hilbert space ${\mathcal H}$ is the space $L^2(M,E)$ of $L^2$-sections of a holonomy
equivariant Hermitian vector bundle $E$, on which an element $k$
of the algebra ${\mathcal A}$ is represented via the
$\ast$-representation $R_E$;
\item The operator $D$ is a first order self-adjoint
transversally elliptic operator with the holonomy invariant
transversal principal symbol such that the operator $D^2$ has the
scalar principal symbol.
\end{enumerate}

Let $S^*{\cA}$ denote the unitary cotangent bundle and $\gamma_t$
the noncommutative geodesic flow associated with $({\cA}, {\cH},
D)$ (see \cite{mpag} for definitions in the non-unital case).
Thus, $S^*{\cA}$ is a $C^*$-algebra and $\gamma_t$ is a
one-parameter group of its automorphisms.

The transversal bicharacteristic flow
$\operatorname{Ad}(\alpha_t)^*$ of the operator $\langle D\rangle$
extends by continuity to a strongly continuous one-parameter group
of automorphisms of the algebra $\bar{S}^0(G_{{\mathcal
F}_N},{\mathcal L}(\pi^*E)\otimes |T{\cG}_N|^{1/2})$, the uniform
closure of $S^0(G_{{\mathcal F}_N},{\mathcal L}(\pi^*E)\otimes
|T{\cG}_N|^{1/2})$ (see \cite{mpag}).

\begin{thm}
\label{noncom:flow} Let $({\mathcal A},{\mathcal H},D)$ be a
spectral triple associated with a compact foliated Riemannian
manifold $(M,{\mathcal F})$ as above. There exists a surjective
homomorphism of involutive algebras
\[
P : S^*{\mathcal A}\rightarrow \bar{S}^0(G_{{\mathcal
F}_N},{\mathcal L}(\pi^*E)\otimes |T{\cG}_N|^{1/2})
\]
such that the following diagram commutes:
\begin{equation}\label{e:cd}
  \begin{CD}
S^*{\mathcal A} @>\gamma_t>> S^*{\mathcal A}\\ @VPVV @VVPV
\\ \bar{S}^0(G_{{\mathcal F}_N},{\mathcal L}(\pi^*E)\otimes|T{\cG}_N|^{1/2})@>\operatorname{Ad}(\alpha_t)^*>>
\bar{S}^0(G_{{\mathcal F}_N},{\mathcal
L}(\pi^*E)\otimes|T{\cG}_N|^{1/2})
  \end{CD}
\end{equation}
\end{thm}

Here the map $P$ is induced by the principal symbol map
$\bar{\sigma}$, and the theorem is a simple consequence of the
results of \cite{mpag} and Theorem~\ref{Egorov}.

\end{document}